\documentclass[11pt,a4paper]{article}
\usepackage{amsfonts,amsmath,amsthm,epsfig,verbatim,amssymb,amscd,eucal,bbm}
\usepackage[latin1]{inputenc}

\theoremstyle{plain}
\newtheorem{thm}{Theorem}[section]
\newtheorem{prop}[thm]{Proposition}
\newtheorem{defn}[thm]{Definition}
\newtheorem{cor}[thm]{Corollary}
\newtheorem{lem}[thm]{Lemma}

\newtheorem{defn-prop}[thm]{Definition-Proposition}

 \newcommand{\cA}{{\cal A}}
 \newcommand{\cB}{\cal{B}}
 \newcommand{\Ce}{\cal{C}}
\newcommand{\cD}{{\cal D}}

\newcommand{\cF}{{\cal F}}
\newcommand{\cI}{{\cal I}}

\newcommand{\cX}{{\cal X}}
\newcommand{\cY}{{\cal Y}}

\newcommand{\mrm}[1]{\text{\rm #1}}

\def\eref#1{(\ref{#1})}
\newcommand{\E}{\mathbb{E}}
\newcommand{\N}{\mathbb{N}}
\newcommand{\p}{\mathbb{P}}

\newcommand{\R}{\mathbb{R}}

\newcommand{\Z}{\mathbb{Z}}

\begin{document}

\noindent

\title{\bf The Single Server Queue and the Storage Model: Large Deviations and Fixed Points}

\author{Moez {\sc DRAIEF}
\thanks{Statistical Laboratory,
Centre for Mathematical Sciences, Wilberforce Road, Cambridge CB3
0WB UK Tel: {\tt +44 1223337945} E-mail: {\tt
M.Draief@statslab.cam.ac.uk} }}

\date{}

\maketitle

\begin{abstract}
\noindent We consider the coupling of a single server queue and a
storage model defined as a Queue/Store model in
\cite{DMO03,DMO04}. We establish that if the input variables,
arrivals at the queue and store, satisfy large deviations
principles and are linked through an {\em exponential tilting}
then the output variables (departures from each system) satisfy
large deviations principles with the same rate function. This
generalizes to the context of large deviations the extension of
Burke's Theorem derived in \cite{DMO03,DMO04}.
\end{abstract}

\medskip

%\renewcommand\abstractname{R\'esum\'e}
%\begin{abstract}
%\noindent Nous établissons un résultat de point fixe dans le cadre
%des grandes déviations pour le modèle File d'attente/Modèle de
%Stockage introduit dans \cite{DMO03,DMO04}. Plus précisément si
%les variables d'entrée du système vérifient des principes de
%grandes déviations avec des fonctions de taux liées par une
%relation d'inclinaison exponentielle alors il en est de même pour
%les variables de sortie dont les fonctions de taux sont les mêmes
%que celles des variables d'entrée.
%\end{abstract}
%
%\medskip


\medskip

%%% ----------------------------------------------------------------------

\section{Introduction}
\noindent A celebrated theorem of Burke \cite{Bu56} asserts that
when a Poisson process $\{a_n,\:n\in\Z\}$ with mean inter-arrival
time $1/\lambda$ is input to a single server queue whose service
times $\{s_n,\:n\in\Z\}$ are i.i.d. exponentials (./M/1/$\infty$
following the Kendall nomenclature) with mean $1/\mu<1/\lambda$,
the equilibrium departure process $\{d_n,\:n\in\Z\}$ from the
queue is also a mean $1/\lambda$ Poisson process. In this sense
the mean $1/\lambda$ Poisson process is a {\em fixed point} for
the ./M/1/$\infty$ operator. In \cite{DMO03,DMO04}, we prove a
generalization of Burke's theorem to the couple of departures and
the sequence of times spent at the very back of the queue
$\{r_n,\:n\in\Z\}$. The sequence $\{r_n,\:n\in\Z\}$ can also be
seen as the sequence of departures from a storage model, where
$s_n$ is the amount of $P$ {\bf S}upplied at slot $n+1$, and $a_n$
is the amount of $P$ {\bf A}sked for at the same slot. More
precisely, we show that under the conditions of Burke's theorem
the couple of {\em output variables} from the Queue/Store model
$\{(d_n,r_n),\:n\in\Z\}$ has the same law as the couple of {\em
input variables} $\{(a_n,s_n),\:n\in\Z\}$.

\medskip

In this paper, we consider the fixed point question at large
deviations scaling. A similar result has been presented by Ganesh
et al. in \cite{GOP03} for discrete-time queues. Assuming that the
input variables satisfy a large deviations principle, we show that
if their laws are linked through an exponential tilting, then the
large deviations principle is preserved by the Queue/Store system
(i.e.\ the output variables satisfy large deviations principles
with the same rate functions as the input variables).

\medskip

This paper is organized as follows: section \ref{se-ldp} is
devoted to giving some background on large deviations that will be
useful for our analysis. The Queue/Store model is presented in
section \ref{se-model}. Section \ref{se-LDPw} focuses on the
workload process. We define the concept of effective bandwidth in
section \ref{se-eband}. Finally, in section \ref{se-ldpout}, we
derive a large deviations principle for the output variables of
the Queue/Store model (\S \ref{sse-opt}) and we give the condition
on the rate functions of the input variables (\S \ref{sse-fix})
which are preserved by the Queue/Store operator.

\section{Large deviations}\label{se-ldp}
Let $(\Omega, \cF,P)$ be a probability space and
$x=\{x_n,n\in\N\}$ be a sequence of i.i.d.\ random variables. We
define the sequence of partial sums $X=\{X_n,n\in \N\}$ by
$$X_n=\sum_{i=1}^nx_i\:.$$ If $\E |x_0|<\infty$, then $X$
satisfies {\em the strong law of large numbers}, i.e.,
$$\lim_{n\to\infty}\frac{X_n}{n}=\E x_0,\quad a.s.$$
We also focus on the fluctuation of a random variable around its
mean. Suppose $\E x_0^2<\infty$ then the {\em central limit
theorem} gives the fluctuations of the scale of $O(1/\sqrt{n})$,
of $X_n/n$ around $\E x_0$. More precisely the sequence
$\sqrt{n}\Bigl(\frac{X_n}{n}-\E x_0\Bigr)$ converges in law a
Gaussian random variable with mean zero and variance the one of
$x_0$. Large deviations theory deals with the fluctuations of the
scale of $O(1)$.

\medskip

\begin{defn}
Let $\cX$ be a real Hausdorff space. A function
$I:\cX\rightarrow\overline{\R}_+=\R_+\cup \{\infty\}$ is a {\em
rate function} if $I$ is lower semi-continuous, i.e. the sets
$\{x:I(x)\leq \alpha\}$ are closed, for all $\alpha\in\R$. In
addition, if these sets are compact then $I$ is a {\em good} rate
function.
\end{defn}

\medskip

\begin{defn}
A sequence $\{x_n,\:n\in\N\}$ in $\cX$ satisfies a {\em large
deviations principle} with the rate function $I
:\cX\rightarrow\overline{\R}_+$, if for each Borel subset $\cB$ of
$\cX$,
\begin{eqnarray*}
-\inf_{x\in {\cB°}} I(x) & \leq &
\liminf_{n\to\infty}\frac{1}{n}\log\p(x_n\in {\cB})\\
&\leq & \limsup_{n\to\infty}\frac{1}{n}\log\p(x_n\in {\cB})\leq
-\inf_{x\in {\overline{\cB}}} {I(x)}\:,
\end{eqnarray*}
where ${\cB°}$ denotes the interior of $\cB$ and  $\bar{\cB}$ its
closure.
\end{defn}

\medskip

Throughout the rest of this section, we are interested in
real-valued random variables.

\medskip

\begin{defn}
The {\em cumulant generating function} of a real-valued random
variable $x$ is given by
$$\Lambda_X(\theta)=\log\E e^{\theta x},\quad\forall \theta
\in\R\:,$$ which can be infinite.
\end{defn}

\medskip

Using Hölder's inequality and Fatou's lemma, we check that
$\Lambda_X$ is convex, lower semi-continuous and for all $\theta$
in the effective domain of $\Lambda_X$ (i.e. $\{\theta\mid
\Lambda_X(\theta)<\infty\}$), we have
\[\Lambda'_X(\theta)=\frac{\E(X_0 e^{\theta
X_0})}{e^{\Lambda_X(\theta)}}\:.\]

\medskip

We define the {\em Legendre transform} of $\Lambda_X$ as
\begin{equation}\label{eq-deflegendre}
I_X(x)=\sup_{\theta\in\R}\{\theta  x-\Lambda_X(\theta)\}\:.
\end{equation}

\medskip

The function $I_X$ is positive, convex, and lower semi-continuous.
If $\Lambda_X$ is finite at the origin then
$$\Lambda_X'(0)=\E X_0,\; \text{ and } I_X(\E X_0)=0\:.$$
Moreover, we easily check that
\begin{eqnarray}\nonumber
I_X(x)=\sup_{\theta\geq 0}\{\theta x-\Lambda_X(\theta)\}&\mrm{ for
} x\geq \E x_0\\ \label{eq-IMean} I_X(x)=\sup_{\theta\leq
0}\{\theta x-\Lambda_X(\theta)\}&\mrm{ for } x\leq \E x_0\:.
\end{eqnarray}

\medskip

We recall the statement of {\em Cramér's theorem} for real-valued
random variables.

\medskip

\begin{thm}\label{th-cramer}
Let $\{x_n,\:n\in\N\}$ be an i.i.d. sequence of real-valued random
variables and let $\Lambda_X$ be its cumulant generating function
which we suppose finite at the origin. Then the sequence of random
variables $\{X_n/n,n\in\N^*\}$, where $X_n=\sum_{i=1}^n x_i$,
satisfies a large deviations principle with rate function $I_X$
defined in \eref{eq-deflegendre}: for each closet subset $F$ of
$\R$,
\begin{equation}\label{eq-ldpup}
\limsup_{n\to\infty}\frac{1}{n}\log\p(S_n\in F)\leq -\inf_{x\in F}
I_X(x)\:.
\end{equation}
and for each open subset $O$ of $\R$
\begin{equation}\label{eq-ldpdn}
\liminf_{n\to\infty}\frac{1}{n}\log\p(S_n\in O)\geq -\inf_{x\in O}
I_X(x)\:.
\end{equation}
\end{thm}

\medskip

The inequalities \eref{eq-ldpup} and \eref{eq-ldpdn} are known as
the upper and lower bounds of the large deviations principle. For
the different extensions of Cramér theorem we refer to
\cite{DeZe93}. We can have a large deviations principle for a
couple of independent random variables.

\medskip

\begin{thm}
Let $\cX$ and $\cY$ be two real Hausdorff spaces, and
$\{x_n,\:n\in\N\}$ and $\{y_n,\:n\in\N\}$ be two sequences
satisfying large deviations principles on $\cX$ and $\cY$, with
rate functions $I_X$ and $I_Y$. Suppose that they are independent,
then the sequence $\{(x_n,y_n),\:n\in\N\}$ satisfies a large
deviations principle on $\cX\times\cY$ with rate function
$$I_{X,Y}(x,y)=I_X(x)+I_Y(y)\:.$$
\end{thm}

\medskip

In this paper, we will establish large deviations principles
relying on indirect methods. Once we have a large deviations
principle for one sequence of random variables, we can
effortlessly obtain large deviations principles for a whole class
of random sequences, namely those obtained via continuous
transformations. We present the {\em contraction principle}, the
tool that enables this.

\medskip

\begin{thm}
If $\{x_n,\:n\in\N\}$ satisfies a large deviations principle on
$\cX$ with rate function $I_X$, and if $f:\cX\rightarrow\cY$ is a
continuous function, then the sequence $\{y_n,\:n\in\N\}$ defined
by $y_n=f(x_n)$ satisfies a large deviations principle on $\cY$
with rate function
$$I_Y(y)=\inf_{x,\:f(x)=y}I_X(x)\:.$$
\end{thm}

\medskip

Sometimes $f$ is  ``almost continuous", i.e. the sequence
$\{y_n,\:n\in\N\}$ is close to a sequence $\{f(x_n),\:n\in\N\}$,
where $f$ is continuous.

\medskip

\begin{defn-prop}
Let $\{x_n,\:n\in \N\}$ and $\{y_n,\:n\in\N\}$ be two sequences on
$\cX$. They are said {\em exponentially equivalent}, if for all
$\epsilon>0$,
$$\limsup_{n\to\infty}\log\p(||x_n-y_n||>\epsilon)=-\infty\:.$$
In this case, if $\{x_n,\:n\in \N\}$ satisfies a large deviation
principle then $\{y_n,\:n\in\N\}$ satisfies a large deviations
principle, and $I_X=I_Y$.
\end{defn-prop}

\medskip

\begin{thm}
If $\{x_n,\:n\in\N\}$ satisfies a large deviations principle on
$\cX$ with rate function $I_X$, and if $\{y_n,\:n\in\N\}$ and
$\{f(x_n),\:n\in\N\}$ are exponentially equivalent, with
$f:\cX\rightarrow\cY$ continuous, then $\{y_n,\:n\in\N\}$
satisfies a large deviations principle on $\cY$ with rate function
$$I_Y(y)=\inf_{x\in\cX,\:f(x)=y}I_X(x)\:.$$
\end{thm}

\medskip

For a thorough presentation of large deviations' results and
applications, we refer to \cite{DeZe93}.

\section{The Model}
\label{se-model}

Let $\cA=\{A_n,\:n \in \Z\}$ be a point process and assume that
$A_0\leq 0
 <A_1$ and $A_n<A_{n+1},\:\forall n\in\Z$.
We define the $\R_+^*$-valued sequence of r.v's $a=\{a_n,\:n \in
\Z\}$ by $a_n =A_{n+1} -A_{n}$. Let $s=\{s_n,\:n \in \Z\}$ be
another $\R_+^*$-valued sequence of r.v's. The sequences $a$ and
$s$ are the {\em input variables} of the model.
\medskip

Define the sequence of r.v.'s $\cD=\{D_n,\:n \in \Z\}$ by
\begin{equation}
\label{Depar.times} D_n = \sup_{k \leq n}\: \Bigl[A_k
+\sum_{i=k}^n s_i \Bigr]\:.
\end{equation}

A priori the $D_n$'s are valued in $\R \cup \{+\infty\}$. Assume
that $a$ and $s$ are such that the $D_n$'s are almost surely
finite. Set $d_n=D_{n+1}-D_n$ and $d=\{d_n,\:n \in \Z\}$. Define
an additional sequence of r.v.'s $r=\{r_n,\:n \in \Z\}$, valued in
$\R_+^*$, by
\begin{eqnarray}\label{reverserv.times}
r_n &=& \min(D_n, A_{n+1})-A_n\:.
\end{eqnarray}
The sequences $d$ and $r$ are the {\em output variables} of the
model. In view of the future analysis, it is convenient to
introduce the sequence $w=\{w_n,\:n \in \Z\}$ of random variables
valued in $\R_+$, defined by
\begin{equation}
\label{Devolindley} w_n = D_n -s_n -A_n = \sup_{k\leq n-1} \Bigl[\
\sum_{i=k} ^{n-1} (s_i-a_i) \Bigr]^+ \:.
\end{equation}
These random variables satisfy the following recursion (Lindley's
equation)
\begin{equation}
\label{Reclindley} w_{n+1} = [w_n + s_n -a_n]^+\:.
\end{equation}
Using the variables $w_n$, we can give alternative definitions of
$D_n$ and $r_n$:
\begin{equation}
\label{backD} \forall l\leq n, \ D_n= \Bigl[w_l + A_l +
\sum_{i=l}^n s_i \Bigr] \vee \max_{l< k \leq n}\: \Bigl[A_k
+\sum_{i=k}^n s_i \Bigr]\:,
\end{equation}
\begin{equation}
\label{back2} r_n= \min\{w_n + s_n, a_n\}=s_n+w_n-w_{n+1}\:.
\end{equation}

We now interpret the variables defined above in two different
contexts: a {\em queueing model} and a {\em storage model}.

\subsection{The single-server queue}\label{sse-queue}

We are concerned with a single server queue where each customer is
characterized by an instant of arrival in the queue and a service
demand. Customers are served upon their arrival in the queue and
in their order of arrival. Since there is a single server, a
customer may have to wait in a buffer before the beginning of its
service. Using Kendall's nomenclature, our model is a
$././1/\infty/$FIFO queue. The customers are numbered by $\Z$
according to their order of arrival in the queue (customer 1 being
the first one to arrive strictly after instant 0). Let $A_n$ be
the instant of {\bf A}rrival of customer $n$ and $s_n$ its {\bf
S}ervice time. Then the variables defined in
\eref{Depar.times}-\eref{Reclindley} have the following
interpretations:
\begin{itemize}
\item $D_n$ is the instant of departure of customer $n$ from the
queue, after completion of its service; $\{d_n,\:n\in\Z\}$ is the
sequence of inter-departure times;
\item $w_n$ is the waiting time of customer $n$ in the buffer
between its arrival and the beginning of its service;
\item $r_n$ is the time
spent by customer $n$ at the {\em very back} of the queue.
\end{itemize}
The variables $\{r_n,\:n\in\Z\}$ are less classical in queueing
theory \cite{GaPr01}.

\subsection{The storage model}\label{sse-storage}
Some product $P$ is supplied, sold and stocked in a store in the
following way: Events occur at integer-valued epochs, called {\em
slots}. At each slot, an amount of $P$ is supplied and an amount
of $P$ is asked for by potential buyers. The rule is to meet all
demands, if possible. The demand of a given slot which is not met
{\em is lost}. The supply of a given slot which is not sold {\em
is not lost} and is stocked for future consideration.

\medskip

Let $s_n$ be the amount of $P$ {\bf S}upplied at slot $n+1$, and
let $a_n$ be the amount of $P$ {\bf A}sked for at the same slot.
In this context, the variables in
\eref{Depar.times}-\eref{Reclindley} can be interpreted as
follows:
\begin{itemize}
\item $w_n$ is the level of the stock at the end of slot $n$. It
evolves according to (\ref{Reclindley});
\item $r_n$ is the demand met at slot $n+1$, see equation (\ref{back2}); it
is the amount of $P$ departing at slot $n+1$;
\end{itemize}
The variables $\{D_n,\:n\in\Z\}$ do not have a natural
interpretation in this model.

\medskip

It is important to remark that while the equations driving the
{\em single server queue} and the {\em storage model} are exactly
the same, the relevant variables are different. The important
variables are the ones corresponding to the departures from the
system. The departures are coded in the variables $
\{d_n,\:n\in\Z\}$ for the single server queue and in the variables
$\{r_n,\:n\in\Z\}$ for the storage model. For a more detailed
discussion of these two models we refer to \cite{DMO04}.

\subsection{Rare events}

We recall that $w_n$ is the waiting time of customer $n$ before it
starts its service (respectively the level of the stock at the end
of slot $n$) and it is given by
\begin{equation}\label{eq-lindleys}
w_{n}=(w_{n-1}+s_{n-1}-a_{n-1})^+=\sup_{m\leq
n}[\sum_{k=m}^{n-1}s_k-a_k]\:.
\end{equation}
We assume that the sequences $a=\{a_n,\:n\in\Z\}$ and
$s=\{s_n,\:n\in\Z\}$ are stationary and ergodic. Under the
stability condition $\E s_0<\E a_0$, we can use the above
expression of $w_n$ to give its asymptotic behaviour using the
large deviations properties of the input variables
$a=\{a_n,\:n\in\Z\}$ and $s=\{s_n,\:n\in\Z\}$. We start with an
example.

\medskip

\paragraph{Example} Suppose that customers arrive in a deterministic fashion, i.e. $a_n=1,\:\forall n$,
and that the sequence $\{s_n,\:n\in\Z\}$ is i.i.d. with
$$\p(s_n=2)=1-\p(s_n=0)=p<\frac{1}{2}\:.$$
The process $\{w_n,\:n\in\Z\}$ is a discrete time birth-and-death
process with stationary distribution
$$\p(w_0 \geq q)=\Bigl(\frac{p}{1-p}\Bigr)^q\:.$$
We get
$$\frac{1}{n}\log \p(w_0 \geq nq)=-\delta q\:,$$
where $\delta =\log\frac{1-p}{p}$.

\medskip

This approximation remains valid under general conditions on the
input variables, with an expression of $\delta$ that depends upon
the rate functions associated to the sequences $a$ and $s$. Indeed
it has been proved \cite{Cha94b,DuOc95,GlWh94} that the stationary
version of $w_n$ satisfies the following property
\begin{equation}\label{eq-overflow}
\lim_{n\to\infty}\frac{1}{n}\log\p(w_0>nq)=-\delta q\:,
\end{equation}
or alternatively
$$\p(w_0/n>q)\asymp e^{-n\delta q}\:.$$
First we give a heuristic proof of this result which will be
asserted rigorously in section \ref{se-LDPw}. We assume the
sequences $\{a_n,\:n\in\Z\}$ and $\{s_n,\:n\in\Z\}$ i.i.d.
mutually independent with
$$\Lambda_A(\theta)=\log\E e^{\theta a_0}\:,\qquad\Lambda_S=\log\E
e^{\theta s_0}\:,$$ finite near the origin. Let $x_n=s_{n}-a_{n}$
for all $n\in\Z$, and $X_k=\sum_{i=1}^{k} x_{-i}$ for all $k\in\N$
then $w_0=\sup_{k\geq 0} X_k$, with $X_0=0$. Using Cramér's
theorem (Theorem \ref{th-cramer}), the sequences
$\{a_n,\:n\in\Z\}$ and $\{s_n,\:n\in\Z\}$ satisfy large deviations
principles on $\R$ with rate function $I_A$ and $I_S$ the Legendre
transforms of $\Lambda_A$ and $\Lambda_S$. Thus the sequence
$\{x_n,\:n\in\Z\}$ satisfies a large deviations principle with
rate function $I_X$, which is the Legendre transform of
$\Lambda_X(\theta)=\Lambda_S(\theta)+\Lambda_A(-\theta)$, i.e.
$\p(\frac{X_n}{n}>x)\asymp e^{-nI_X(x)}$. Moreover,
\begin{eqnarray*}
\p(w_0\geq n q)=\p(\sup_{k\geq 0} X_k\geq n q)
&=&\p(\cup_{k\geq 0} \{X_k\geq q n\})\\
&\leq&\sum_{k=1}^{\infty}\p(X_k\geq q n)\:.
\end{eqnarray*}
Since $\p(X_k\geq n q)=\p(\frac{X_k}{k}\geq \frac{n q}{k})\asymp
e^{-kI_X(nq/k)}$, we check that
\begin{equation}\label{eq-LDPw}
\p(w_0\geq
nq)\asymp\sum_{k=1}^{\infty}e^{-nq\frac{I_X(nq/k)}{nq/k}}\:.
\end{equation}
We conclude using the principle of the largest term, i.e.
$$\sum_{k=1}^{\infty}e^{-nq\frac{I_X(nq/k)}{nq/k}}\asymp
e^{-qn\delta}$$ where $\delta=\inf_{x>0}\frac{I_X(x)}{x}$ and we
get \eref{eq-overflow}. The principle of the largest term
translates the fact that rare events occur in the most likely way.
Indeed, the dominant term in \eref{eq-LDPw}, which gives the
explosion of the waiting times in the queueing system (or the
overflow of the stock in the storage model), happens following the
most probable scenario.

\section{Large deviations for the workload process}\label{se-LDPw}
First, we give the assumptions under which the approximation
\eref{eq-overflow} is fulfilled. For more general assumptions we
refer to the paper by Ganesh et al. \cite{GOW03}.

\medskip

\paragraph{Assumptions:}
\begin{itemize}
\item [$(i)$]  The sequences $\{a_n,\:n\in\Z\}$ and
$\{s_n,\:n\in\Z\}$ are i.i.d., mutually independent. Their
cumulant generating functions are given by
\[
\Lambda_A(\theta)=\log \E e^{\theta a_0},\qquad
\Lambda_S(\theta)=\log \E e^{\theta s_0}
\]
We assume that both $\Lambda_A$ and $\Lambda_S$ are differentiable
near the origin.
\item [$(ii)$] The stability condition $$1/\mu=\Lambda'_S(0)=\E s_0<\E
a_0=\Lambda'_A(0)=1/\lambda$$
is satisfied.
\end{itemize}

\medskip

Under $(i)$ and $(ii)$, we have

\medskip

\begin{prop}\label{PGDw}
The sequence $\{w_0/n,\:n\in\N^*\}$ satisfies a large deviations
principle with good rate function $I_W(q)=\delta q$, i.e.
$\lim_{n\to\infty}\frac{1}{n}\log\p(w_0/n\geq q)=-\delta q$, with
\begin{equation}\label{eq-delta}
\delta=\inf_{0<a<s}\frac{I_{A,S}(a,s)}{s-a}=\sup\{\theta:\:\Lambda_S(\theta)+\Lambda_A(-\theta)\leq
0\}\:.
\end{equation}
\end{prop}

\medskip

\begin{proof}
First we prove \eref{eq-delta}. Recall that $x_n=s_n-a_n$,
$X_n=\sum_{i=1}^{n}x_{-i}$, and
$$w_0=\sup_{k\geq 0} X_k\:.$$
By the independence assumption, we have
$$\Lambda_X(\theta)=
\Lambda_S(\theta)+\Lambda_A(-\theta)\:.$$ By direct application of
the contraction principle, we check that
$$I_X(x)=\inf_{y> x}\{I_S(y)+I_A(y-x)\}\:.$$ Let $\theta \leq
\inf_{x\geq 0} I_X(x)/x$ then
\begin{eqnarray*}
\theta \leq \inf_{x\geq 0} I_X(x)/x &\Leftrightarrow&
\theta \leq I_X(x)/x, \forall x\geq 0\\
&\Leftrightarrow&
\theta x -I_X(x)\leq 0, \forall x\geq 0\\
&\Leftrightarrow&
\sup_{x\geq 0} \{\theta x -I_X(x)\} \leq 0\\
&\Leftrightarrow& \Lambda_X(\theta) \leq 0\:.
\end{eqnarray*}
The last equivalence is due to the fact that $\E x_0<0$ (stability
condition) and the equation \eref{eq-IMean}. We proved that
\[
\inf_{0<a<s}\frac{I_{A,S}(a,s)}{s-a}=\inf_{x\geq 0}
I_X(x)/x=\sup\{\theta:\:\Lambda_S(\theta)+\Lambda_A(-\theta)\leq
0\}\:.
\]

\medskip

$\bullet$ {\em Lower bound:} For $q>0$, we have $\p(w_0\geq q)\geq
P(X_k\geq q)$. Notice that for $p\geq\frac{q}{\lceil q/p\rceil}$,
we have
$$\p(w_0\geq q)\geq\p\Bigl(\frac{1}{\lceil q/p \rceil} X_{\lceil q/p
\rceil}\geq p\Bigr)\:.$$ Since $\frac{1}{q}\geq
\frac{1}{p}\frac{1}{\lceil q/p\rceil}$, we get
\begin{eqnarray*}
\liminf_{q\to\infty}\frac{1}{q}\log\p(w_0\geq
q)&\geq&\frac{1}{p}\liminf_{q\to\infty}\frac{1}{\lceil q/p
\rceil}\log\p\Bigl(\frac{1}{\lceil q/p \rceil} X_{\lceil q/p
\rceil}\geq p\Bigr)\\
&=&\frac{1}{p}\liminf_{n\to\infty}\frac{1}{n}\log
\p\Bigl(\frac{1}{n} X_{n}\geq p\Bigr)\geq -\frac{1}{p}I_X(p)\:,
\end{eqnarray*}
for all $p>0$. We conclude that
\begin{eqnarray*}
\liminf_{q\to\infty}\frac{1}{q}\log\p(w_0\geq q)\geq
-\inf_{x>0}\frac{1}{x}I_X(x)&=&-\inf_{x>0}\frac{1}{x}\inf_{y>
x}\{I_S(y)+I_A(y-x)\}\\
&=&-\inf_{0<a<s}\frac{I_S(s)+I_A(a)}{s-a}\:.
\end{eqnarray*}

\medskip

$\bullet$ {\em Upper bound:} Recall that
$\Lambda_X(\theta)=\Lambda_S(\theta)+\Lambda_A(-\theta)$. By the
stability condition, we have
$$\Lambda_X'(0)=\log\frac{\E(s_1)}{\E(a_1)}<0\:.$$ Moreover,
 $\Lambda_X$ is differentiable near
$0$ with $\Lambda_X(0)=0$. Thus, there is a constant $\Theta>0$
such that $\Lambda_X(\Theta)<0$ and $\E e^{\Theta X_n}<\infty$.
Applying Chernoff's bound, we get, for $n\in\N$
$$\p(X_n\geq q)\leq e^{-\Theta q}\E e^{\Theta X_n}\:.$$
This leads to
$$\limsup_{q\to\infty}\frac{1}{q}\log\p(X_n\geq q)\leq -\Theta\:.$$
For $N\in\N$, we check that
$$\p\Bigl(\max_{0\leq n\leq N}X_n\geq q\Bigr)\leq N \max_{0\leq n\leq
N}\p(X_n\geq q)\:.$$ Allowing $q$ go to infinity,
$$\limsup_{q\to\infty}\frac{1}{q}\log\p\Bigl(\max_{0\leq n\leq N}X_n\geq q\Bigr)\leq
\max_{0\leq n\leq N}\limsup_{q\to\infty}\log\frac{1}{q}\p(X_n\geq
q)\leq -\Theta\:.$$ We need now to have a bound for large values
of $n$. Applying the union bound, we have
\begin{equation}\label{eq-unionbd}
\p(\sup_{n>N}X_n\geq q)\leq \sum_{n>N}\p(X_n\geq q)\leq e^{-\Theta
q}\sum_{n>N}\E  e^{\Theta X_n}\:.
\end{equation}\label{eq-upbd}
Since $\frac{1}{n}\log\E e^{\Theta X_n}=\Lambda_X(\Theta)<0$,
there is $0<\epsilon<-\Lambda_X(\Theta)$ et $N_{\Theta}\in \N$
such that
\begin{equation}
\forall n>N_{\Theta},\:.\frac{1}{n}\log\E e^{\Theta X_n}\leq
-\epsilon\:.
\end{equation}
Combining \eref{eq-unionbd} and \eref{eq-upbd}, we obtain
$$\p\Bigl(\sup_{n>N_{\Theta}}X_n>q\Bigr)\leq e^{-\Theta
q}\sum_{n>N_{\Theta}}e^{-n\epsilon}<\frac{e^{-\Theta
q}}{1-e^{-\epsilon}}e^{-(N_{\Theta} +1)\epsilon}\:.$$ To sum up,
we derived the following bound
$$\limsup_{q\to\infty}\frac{1}{q}\log\p\Bigl(\sup_{n>N_{\Theta}}X_n\geq\Bigr)\leq-\Theta\:.$$
Therefore, for
 $\Theta>0$ such that $\Lambda_X(\Theta)<0$, we have
$$\limsup_{q\to\infty}\frac{1}{q}\log \p\Bigl(\sup_{n\geq
0}X_n\geq q\Bigr)\leq -\Theta\:.$$ We conclude using
\eref{eq-delta}.

\end{proof}

This result can be interpreted naively in terms of the following
approximation
$$\p(\sup_n X_n\geq q)\approx \sup_n\p( X_n\geq q)\:.$$
This is nothing more than a consequence of the fact that rare
events occur in the most likely way. Moreover, using this
approximation one can predict the frequency at which the waiting
time (or the level of stock) overflows a given threshold. Indeed
if \eref{eq-overflow} is fulfilled, then we represent the graph
(in logarithmic squale) of the frequency of exceeding a level $q$.
The above approximation (Proposition \ref{PGDw}) inspired numerous
applications (described by Courcoubetis et al. \cite{CKRWW95}) of
large deviations to the analysis of the statistics of networks of
queues.

\section{Effective bandwidth}\label{se-eband}
The {\em bandwidth} of a traffic flow has become a classical
feature in communication networks literature. We are interested in
the notion of {\em effective bandwidth} introduced by Kelly in
\cite{Kel91}, which gives an analytical way of describing the
properties of a stochastic flow by means of rare events occurring
in the network through which it passes.

\medskip

In practice, the queues and stores we are interested in have
finite capacities, i.e.\ the queue (respectively the store)
rejects customers (respectively stock) each time the waiting time
(respectively the amount of stock) overpasses a given threshold.
In the queueing context, we assume $\{s_n,\:n\in \Z\}$ given and
we seek the (deterministic) minimal value $a_p$ of inter-arrival
times such that the probability that the waiting times exceeds the
threshold is less than a fixed value $p$. For the storage model,
we assume $\{a_n,\:n\in \Z\}$ given and we seek the maximal value
$s_p$ of product arriving at the store at each slot of time, such
that the store rejects it with a probability less then $p$. More
precisely, we want to identify the arrivals to both models such
that
\begin{equation}\label{eq-overload}
\p(w_0\geq q)\leq p\:,
\end{equation}
for given values of $q$ and $p$.

\paragraph{Queue} Let $a_n=a,\:\forall n\in\Z$, i.e.
$\Lambda_A(\theta)=a \theta$. If $s$ is i.i.d., then by
Proposition \ref{PGDw},
$$\p(w_0\geq q)\leq e^{-\delta(a)q}\:$$
with
$$\delta(a)=\sup\{\theta:\: \theta\geq 0 \mrm{ and }\
\Lambda_S(\theta) \leq a\theta\}\:.$$ The minimal inter-arrival
time for the inequality \eref{eq-overload} to be fulfilled is
given by $a_p$ such that
$$a_p =\inf\{a:\:a\geq 0 \mrm{ et } e^{-\delta(a) q}\leq p\}\:.$$

\medskip

Let $\theta_p=-\frac{\log p}{q}$, then $\delta(a_p)=\theta_p$ and
$$a_p=\frac{\Lambda_S(\theta_p)}{\theta_p}\:.$$
The variable $a_p$ is the {\em effective bandwidth} of the queue.
For more details we refer to Kelly's review \cite{Kel96}.

\medskip

\paragraph{Store} Let $s_n=s,\:\forall n\in\Z$, then
$\Lambda_S(\theta)=\theta s$ and
$$\delta(s)=-\inf\{\theta<0:\ \ \Lambda_A(\theta)\leq s\theta\}\:.$$
the inequality \eref{eq-overload} is satisfied if no more than
$s_p$ amount of stock arrives to the store at each slot with
$$s_p =\sup\{s\geq 0: e^{-\delta(s) q}\leq p\}\:.$$

\medskip

The quantity $s_p$ is the {\em effective bandwidth} of the store,
with
$$s_p=\frac{\Lambda_A(\theta_p)}{\theta_p}\:.$$

\medskip

We generalize this notion to non-deterministic arrivals by
introducing the functions $\alpha_A(\theta)$ and
$\alpha_S(\theta)$ given by
$$\alpha_A(\theta)=\frac{\Lambda_S(\theta)}{\theta},\qquad\alpha_S(\theta)=\frac{\Lambda_A(\theta)}{\theta}\:,$$
representing the effective bandwidths of the queue and the store,
respectively.

\section{Principal of large deviations for the output
variables}\label{se-ldpout} We proved, in section \ref{se-LDPw}, a
large deviations principle for the workload process using
classical techniques. It is generally hard to apply these
techniques to derive large deviations principles for the output
variables. In this section, we will apply the contraction
principle to this end. Following the outline of the proof of the
existence of fixed points for discrete-time queues in
\cite{GOP03}, we use theoretical results on large deviations for
continuous-time processes.

\medskip

First, we go back to the workload process to illustrate this
method. For $n\in\N^*$, we define
$$A_n=\sum_{i=1}^{n}a_{-i},\qquad S_n=\sum_{i=1}^{n} s_{-i}\:.$$
It is obvious that $w_0=\sup_{k\geq 0}(S_k-A_k)$.

\medskip

To apply the contraction principle it is crucial to define an
adapted topology under which the workload process is obtained
through a continuous mapping of the input variables. Therefore, we
define the polygonal approximation in $n$ of a given sequence
$\{X_k,\:k\in\N^*\}$ by
$$\tilde{X}_n(t)=\frac{1}{n}X_{\lfloor nt\rfloor}+(t-\frac{\lfloor
nt\rfloor}{n})(X_{\lfloor nt\rfloor+1}-X_{\lfloor
nt\rfloor}),\qquad  t \geq 0\:.$$

Let $\tilde{A}_n$ and $\tilde{S}_n$ the polygonal approximations
in $n$ of the sequences $\{A_k,\:k\in\N\}$ and $\{S_k,\:k\in\N\}$.
Notice that for $n\in\N^*$,
\begin{eqnarray*}
\frac{w_0}{n}=\frac{1}{n}\sup_{k\in\N}(S_k-A_k)
&=&\frac{1}{n}\sup_{t\geq
0}(S_{\lfloor nt \rfloor}-A_{\lfloor nt\rfloor})\\
&=&\sup_{t>0}(\tilde{S}_n(t)-\tilde{A}_n(t))\:.
\end{eqnarray*}

\medskip

A function $x$ is absolutely continuous (on $\R$) if for every
$\epsilon>0$, there is $\eta>0$  such that for any finitely many
disjoint open intervals $(a_i,b_i),\: i=1\dots n$, satisfying
$\sum_{i=1}^n (b_i-a_i)\leq \eta$, we have $\sum_{i=1}^n
(f(b_i)-f(a_i))\leq \epsilon\:.$ If $x$ is absolutely continuous,
then its derivative $x'$ exists almost everywhere and we can write
$$x(v)-x(u)=\int_u^v x'(t) dt\:.$$

\medskip

For $\mu>0$, we  define $\Ce_{\mu}$ (respectively $\cA_{\mu}$) the
set of continuous functions (respectively absolutely continuous)
$x:\R_+\rightarrow \R$ with $x(0)=0$ and
$$\lim_{t\to+\infty}\frac{x(t)}{t+1}=1/\mu<\infty\:,$$ equipped
with the norm
\begin{equation}\label{norm}
||x||=\sup_{t\in\R_+}|\frac{x(t)}{t+1}|.
\end{equation}
We now focus on large deviations on processes in continuous  time
in $\Ce_{\mu}$.

\medskip

\begin{defn}
A sequence of processes $\{X_n,\:n\in\N\}$, where
$X_n\in\Ce_{\mu}$ satisfies a {\em functional large deviations
principle with linear geodesics}, with instantaneous rate function
$I$, if
\begin{itemize}
\item[$(i)$] the function $I$ is a rate function, with $I(1/\mu)=0$,
\item[$(ii)$] the sequence $\{X_n,\:n\in\N\}$ satisfies a large deviations principle $\Ce_{\mu}$
with rate function
\[
\cI_\mathrm{X}(\phi) = \begin{cases} \int_0^{+\infty}I(\phi'(t))dt  & \text{if } \phi \in \cA_{\mu}\\
                        +\infty & \text{otherwise. }

\end{cases}
\]
\end{itemize}
\end{defn}

\medskip

Dembo and Zajic explored extensions \cite{DeZa95} of large
deviations principles to processes in continuous time under the
topology of uniform convergence on compact intervals. However,
this topology is not appropriate in the context of queueing
systems. We introduce a coarser topology corresponding to the norm
\eref{norm}. The following large deviations principle is due to
Ganesh and O'Connell \cite{GaOc99}. Notice that if
$x=\{x_n,\:n\in\N\}$ an i.i.d. sequence with mean $1/\mu$, then
for $n\in\N^*$, $\tilde{X}_n\in\Ce_{\mu}$, where $\tilde{X}_n$ is
the polygonal approximation in $n$ of the sequence
$X=\{X_n,n\in\N\}$ with $X_n=\sum_{i=1}^n x_i$

\medskip

\begin{thm}\cite{GaOc99}
Let $x=\{x_n,\:n\in\N\}$ be a sequence of real-valued random
variables and $\Lambda_X$ its cumulant generating function, which
we suppose differentiable near the origin. Then the sequence of
polygonal approximations $\{\tilde{X}_n,\:n\in\N^*\}$ of
$X=\{X_n,\:n\in\N\}$, where $X_n=\sum_{i=1}^n x_i$, satisfies a
functional large deviations principle with linear geodesics on
$\Ce_{\mu}$ equipped with the norm \eref{norm}, with mean
$1/\mu=\Lambda_X'(0)$ and instantaneous rate function $I_X$, i.e.
\[
\cI_\mathrm{X}(\phi) = \begin{cases} \int_0^{+\infty}I_X(\phi'(t))dt  & \text{if } \phi \in \cA_{\mu}\\
                        +\infty & \text{otherwise. }

\end{cases}
\]
\end{thm}

\medskip

Under $(i)$, the processes $\tilde{\mathrm{A}}_n$ and
$\tilde{\mathrm{S}}_n$ satisfy functional large deviations
principles with linear geodesics and instantaneous rate functions
$\cI_\mathrm{A}$ and $\cI_\mathrm{S}$ defined by
\[
\cI_\mathrm{A}(\phi) = \begin{cases} \int_0^{+\infty}I_A(\phi'(t))dt  & \text{if } \phi \in \cA_{\lambda}\\
                        +\infty & \text{otherwise , }

\end{cases}
\]
with $I_A(x)=\sup_{\theta\in \R}\{\theta x-\Lambda_A(\theta)\}$,
the Legendre transform of $\Lambda_A$. We define $I_S$ and
$\cI_\mathrm{S}$ in the same way.

\medskip

We check (see \cite{GOW03}) that $f:\Ce_{\mu}\times
\Ce_{\lambda}\rightarrow \R_+$ defined by
$$f(\phi,\psi)=\sup_{t>0}(\psi(t)-\phi(t))\:,$$ is continuous on
$\Ce_{\mu}\times \Ce_{\lambda}$ for the topology induced by the
norm \eref{norm}. Indeed  $f$ is not continuous under the topology
of uniform convergence ( \cite[Example 5.2]{GOW03}).

\medskip

Applying the contraction principle, we check that
$\{w_1/n,n\in\Z\}$ satisfies a large deviations principle on
$\R_+$ with rate function
\begin{eqnarray*}
I_W(q)&=&\inf\{\int_0^{\infty}I_{A,S}(\phi'(s))\:ds\mid
(\phi_1,\phi_2)\in \cA_{\lambda}\times
\cA_{\mu},\:f(\phi_1,\phi_2)=q\}\\
&=&\inf\{\int_0^{\infty}I_{A,S}(\phi'(s))\:ds\mid
(\phi_1,\phi_2)\in \cA_{\lambda}\times
\cA_{\mu},\:\sup_{t>0}(\phi_2(t)-\phi_1(t))=q\}\:,
\end{eqnarray*}
where $I_{A,S}(\phi'(s))=I_{A}(\phi_1'(s))+I_{S}(\phi_2'(s))$.
Combining this with Proposition \ref{PGDw}, we have

\medskip

\begin{cor}\label{le-OptPGDw}
The solution of the following optimization problem
\[\begin{cases}
\text{Minimizing} \int_0^{\infty}I_{A,S}(\phi'(s))\:ds \mrm{ pour
} \phi=(\phi_1,\phi_2)\in
\cA_{\lambda}\times \cA_{\mu}\\
\text{subject to } \sup_{t>0}(\phi_2(t)-\phi_1(t))=q
\end{cases}
\]

is given by $I_W(q)=\delta q$ where $\delta$ is defined in
\eref{eq-delta}.
\end{cor}

\medskip

In the following paragraph, we make use of the same techniques to
give large deviations principles for the output variables
\begin{equation}\label{eq-output}
d_n=a_n+w_{n+1}-w_n+s_{n+1}-s_n, \qquad r_n=s_n-w_{n+1}+w_n\:.
\end{equation}
We define $D=\{D_n,\:n\in\N\}$ and $R=\{R_n,\:n\in\N\}$ by
$$D_n=\sum_{i=1}^n d_{-i},\qquad R_n=\sum_{i=1}^n r_{-i}\:.$$
Let $\tilde{D}_n$ and $\tilde{R}_n$ be their polygonal
approximations in $n$. Under $(i)$, $(ii)$ and by Loynes' scheme
\cite{Loy62}, we check that
$$(\tilde{D}_n,\tilde{R}_n)\in\Ce_{\lambda}\times\Ce_{\mu}\:.$$

    \subsection{Optimization problem}\label{sse-opt}
To make the proof as clear as possible, we introduce a new
variable
$$b_n=a_n+w_{n+1}-w_n\:,$$
which stands, in a queueing context, for the time between
instances of beginning of services for customers $n$ and $n+1$.
Let $B=\{B_n,\:n\in\N\}$ defined by $B_n=\sum_{i=1}^nb_{-i}$ and
$\tilde{B}_n$ its polygonal approximation in $n$ then
$\tilde{B}_n\in\Ce_{\lambda}$. Under $(i)$, we have

\medskip

\begin{lem}
for $0\leq t\leq 1$ fixed, the sequences
$\{\tilde{D}_n(t)\:,n\in\N^*\}$ and
$\{\tilde{B}_n(t)\:,n\in\N^*\}$ are exponentially equivalent.
\end{lem}

\medskip

\begin{proof}
Let $\gamma>0$,
\begin{eqnarray*}
\p(\sup_{0\leq t \leq
1}|\tilde{D}_n(t)-\tilde{B}_n(t)|>\gamma)&\leq
&\p(\sup_{0\leq k\leq n}|s_{-k}-s_{0}|> n \gamma)\\
 &\leq& \sum_{k=0}^n\p(|s_{-k}-s_{0}|>n \gamma)\\
&\leq& e^{-\gamma\nu n}\sum_{k=0}^n \E e^{\gamma|s_{-k}-s_{0}|}
 \leq (n+1) e^{-\gamma\nu n} e^{2\Lambda_S(\gamma)}\:.
\end{eqnarray*}
Since $\lim_{n\to\infty} \log ((n+1) e^{-\gamma\nu n})=-\infty$
and thanks to assumption $(i)$ ($\Lambda_S(\gamma)<\infty$ for
$\gamma$ close to $0$), we are done.
\end{proof}
Thus if one of the sequences $\{\tilde{D}_n(t)\:,n\in\N^*\}$ and
$\{\tilde{B}_n(t)\:,n\in\N^*\}$ satisfies a large deviations
principle then the other does too and they  both have the same
rate function. Since the expression of $\tilde{B}_n$ (in terms of
the input variables) is easier than the one of  $\tilde{D}_n$, we
concentrate on $\tilde{B}_n$ and $\tilde{R}_n$. In fact,

\medskip

\begin{lem}\label{le-Phi}
The sequences $\tilde{B}_n(t)$ and $\tilde{R}_n(t)$ satisfy
\begin{eqnarray}\nonumber
\tilde{B}_n(t)&=&\tilde{A}_n(t)+\sup_{s>0}\{\tilde{S}_n(s)-\tilde{A}_n(s)\}
-\sup_{s>t}\{(\tilde{S}_n(s)-\tilde{S}_n(t))-(\tilde{A}_n(s)-\tilde{A}_n(t))\}\\\label{eq-foncinout}
\nonumber
\tilde{R}_n(t)&=&\tilde{S}_n(t)-\sup_{s>0}\{\tilde{S}_n(s)-\tilde{A}_n(s)\}
+\sup_{s>t}\{(\tilde{S}_n(s)-\tilde{S}_n(t))-(\tilde{A}_n(s)-\tilde{A}_n(t))\}
\end{eqnarray}
\end{lem}

\medskip

\begin{proof}
We check the result for $\tilde{B}_n(t)$, the proof goes along the
same arguments for $\tilde{R}_n(t)$. Since $\tilde{B}_n(t)$ is a
linear interpolation of $B_{\lfloor nt \rfloor}/n$, we will prove
the result for $B_{\lfloor nt \rfloor}$. First, we check that
$$B_{\lfloor nt \rfloor}=A_{\lfloor nt \rfloor}+w_0-w_{-\lfloor nt
\rfloor}\:.$$ However, we have
$w_0=\sup_{s>0}\{\tilde{S}_n(s)-\tilde{A}_n(s)\}$ then
$$w_{-\lfloor nt
\rfloor}=\sup_{s>t}\{(\tilde{S}_n(s)-\tilde{S}_n(t))-(\tilde{A}_n(s)-\tilde{A}_n(t))\}\:.$$
\end{proof}

\medskip

\begin{thm}\label{flle-action}
Under assumptions $(i)$ and $(ii)$, the sequence
$\{(\frac{D_n}{n},\frac{R_n}{n},\frac{w_0}{n}),\:n\in\N^*\}$
satisfies a large deviations principle on $\R_+^3$ with rate
function
\begin{equation}\label{eq-Jrate}
J(x_1,x_2,w)=\inf\{\delta (w-x_1+x_2)+I_{A,S}(x_2,x_1);\inf_C
g(q,\tau,v_1,v_2)\}\:,
\end{equation}
where
\[
g(q,\tau,v_1,v_2)=\tau
I_{A,S}\Bigl(\frac{x_2-v_2}{\tau},\frac{x_2-v_2+w}{\tau}\Bigr)
+(1-\tau
)I_{A,S}\Bigl(\frac{v_1}{1-\tau},\frac{v_2-q}{1-\tau}\Bigr)+\delta
q\:,
\]
with $\delta=\inf_{0<a<s}\frac{I_{A}(a)+I_S(s)}{s-a}$ and
$C=\{(q,\tau,v_1,v_2)\:\mid\:\:q\geq 0,\:0\leq \tau\leq 1,\:
x_2-v_2+w=x_1-v_1+q\}$.
\end{thm}

\medskip

\begin{proof}
By Loynes \cite{Loy62} and lemma \ref{le-Phi}, we check that
$(\tilde{B}_n,\tilde{R}_n,\frac{w_0}{n})$
=$\Phi(\tilde{A}_n,\tilde{S}_n)$ with
$$\Phi:\Ce_{\lambda}\times\Ce_{\mu}
\rightarrow\Ce_{\lambda}\times\Ce_{\mu}\times\R_+\:.$$

Let $\Phi=(\Phi_1,\Phi_2,\Phi_3)$ and
$\phi=(\phi_1,\phi_2)\in\Ce_{\lambda}\times\Ce_{\mu}$, we have
\begin{eqnarray}\nonumber
\Phi_1(\phi)(t)&=&\phi_1(t)+\sup_{s>0}\{\phi_2(s)-\phi_1(s)\}
-\sup_{s>t}\{(\phi_2(s)-\phi_2(t))-(\phi_1(s)-\phi_1(t))\}\\
\nonumber
\Phi_2(\phi)(t)&=&\phi_2(t)-\sup_{s>0}\{\phi_2(s)-\phi_1(s)\}
+\sup_{s>t}\{(\phi_2(s)-\phi_2(t))-(\phi_1(s)-\phi_1(t))\}\\
\label{eq-outin}
\Phi_3(\phi)(t)&=&\sup_{s>0}\{\phi_2(s)-\phi_1(s)\}\:.
\end{eqnarray}
Applying the contraction principle, the sequence
$\{(\tilde{B}_n,\tilde{R}_n,w_0/n),\:n\in\N^*\}$ satisfies a large
deviations principle with rate function
\begin{equation*}
J(\psi_1,\psi_2,w)=\inf\{\int_0^{\infty}I_{A,S}(\phi_1'(s),\phi_2'(s))ds\:\mid\Phi(\phi_1,\phi_2)=(\psi_1,\psi_2,w)\}\:.
\end{equation*}
This variational problem is generally hard to solve, we restrict
ourselves to $t=1$. By the law of large numbers,
$\{(\tilde{B}_n(1),\tilde{R}_n(1)),\:n\in\N^*\}=\{(\frac{B_n}{n},\frac{R_n}{n}),\:n\in\N^*\}$
converges to the mean values of the output variables, whereas the
large deviations principle gives the fluctuations around these
values. The sequence
$\{(\frac{B_n}{n},\frac{R_n}{n},w_0/n),\:n\in\N^*\}$ satisfies a
large deviations principle with rate function
\begin{equation}\label{PGD}
J(x_1,x_2,w)=\inf\{\int_0^{\infty}I_{A,S}(\phi_1'(s),\phi_2'(s))ds\:\mid
[\Phi(\phi_1,\phi_2)](1)=(x_1,x_2,w)\}\:.
\end{equation}
We introduce the variables $q_0$ and $q_1$ given by
\begin{eqnarray*}
q_0&=&\sup_{s>0}\{\phi_2(s)-\phi_1(s)\}\\
q_1&=&\sup_{s>1}\{(\phi_2(s)-\phi_2(1))-(\phi_1(s)-\phi_1(1))\}\:.
\end{eqnarray*} We rewrite
\eref{eq-outin} as follows
\[
x_1=\phi_1(1)+q_0 -q_1,\quad  x_2=\phi_2(1)-q_0 +q_1,\quad w=q_0
\:.
\]
To solve this variational problem we condition on the value of
$q_1$ to be equal to $q$. By the stability condition
$\lambda<\mu$, we check that $q$ is finite. Then $J(x_1,x_2,w)$ is
the solution of the following optimization problem
\begin{eqnarray*}
\mrm{Minimize}&&\int_0^{1}I_{A,S}(\phi'(s))ds+\int_1^{\infty}I_{A,S}(\phi'(s))ds\quad
\mrm{over } \phi\in \cA_{\lambda}\times \cA_{\mu},\: q\geq 0\\
\mrm{subject to}&& q_1=q,\quad x_1=\phi_1(1)+q_0 -q_1,\quad
x_2=\phi_2(1)-q_0 +q_1,\: w=q_0\:,
\end{eqnarray*}
where $\phi=(\phi_1(s),\phi_2(s))$ and
$I_{A,S}(\phi'(s))=I_{A}(\phi_1'(s))+I_{S}(\phi_2'(s))$.

\medskip

Using the auxiliary variable $q$ allows the decomposition of the
initial optimization problem into two minimization problems on
$s\in [0,1]$ and $s>1$ coupled through $q$. More precisely, we
suppose the variables $x_1$, $x_2$ and $w=q_0$ given, if we fix
$q=q_1$, then  $\phi_1(1)$ and $\phi_2(1)$ are fixed. For $s>1$,
the only constraint on $\phi_1$ and $\phi_2$ is $q_1=q$.
Consequently, the first minimization consists of
\[
\mrm{Minimizing }
\int_1^{\infty}I_{A,S}(\phi_1'(s),\phi_2'(s))ds,\quad \mrm{
subject to }q_1=q\:.
\]
The second minimization consists of
\begin{eqnarray*}
\mrm{Minimizing}&& \int_0^{1}I_{A,S}(\phi_1'(s),\phi_2'(s))ds\\
\mrm{ subject to }&& x_1=\phi_1(1)+q_0 -q_1,\quad
x_2=\phi_2(1)-q_0 +q_1,\quad w=q_0\:,
\end{eqnarray*}
knowing $q_1=q$. In both cases the minimization is for $\phi\in
\cA_{\lambda}\times\cA_{\mu}$.

\medskip

For the first problem, let $\psi\in \cA_{\lambda}\times\cA_{\mu}$
with $\psi_i(s-1)=\phi_i(s)-\phi_i(1)$, $s \geq 1,\:i=1,2$. By
lemma \ref{le-OptPGDw}, this optimization gives the rate function
associated with $w_0$ given in Proposition \ref{PGDw} and the
minimum is $\delta q$ where
$\delta=\inf_{a<s}\frac{I_{A,S}(a,s)}{s-a}$.

\medskip

Hence, to determine $J(x_1,x_2,w)$ it remains to
\begin{equation}\label{Min}
\mrm{Minimize }\int_0^1I_A(\phi_1'(s))+I_S(\phi_2'(s))ds+\delta
q\:,
\end{equation}
for $(\phi_1,\phi_2) \in \cA_{\lambda}\times\cA_{\mu}$ subject to
\[
\begin{cases} q&\geq 0,\ \
w= \sup_{s\geq 0}\{\phi_2(s)-\phi_1(s)\} \\
    x_1&=\phi_1(1)+w-q \\
    x_2&=\phi_2(1)-w+q. \\
\end{cases}
\]

\begin{itemize}
\item[$\bullet$] Suppose $\sup_{s\geq
0}\{\phi_2(s)-\phi_1(s)\}=\sup_{s\geq 1}\{\phi_2(s)-\phi_1(s)\}$ ,
then $w=q+\phi_2(1)-\phi_2(1)$, i.e.
$$x_1=\phi_2(1),\quad x_2=\phi_1(1),\quad q=w-x_1+x_2\:.$$ The optimisation problem \eref{Min} is fulfilled by linear
functions $\phi_1$ and $\phi_2$ and we have
\begin{equation*}
\delta (w-x_1+x_2) +I_A(x_2)+I_S(x_1).
\end{equation*}
\item[$\bullet$] If $\phi_2(s)-\phi_1(s)$ reaches its maximum for $\tau\in
[0,1]$ then
\begin{eqnarray*}
x_1&=&-q+\phi_1(1)-\phi_1(\tau)+\phi_2(\tau)\\
x_2&=&q+\phi_2(1)-\phi_2(\tau)+\phi_1(\tau)\:.
\end{eqnarray*}
We define $y_i=\frac{\phi_i(\tau)}{\tau}$ ,
$z_i=\frac{\phi_i(1)-\phi_i(\tau)}{1-\tau},\:i=1,2$ and
 $\chi\in \cA_{\lambda}\times\cA_{\mu}$,
\[
\chi_i(t) = \begin{cases}y_it  & \text{if } t \in ]0,\tau] \\
                        y_i\tau+z_i(t-\tau) & \text{if } t \in ]\tau,1]
\end{cases}\:.
\]
Using the convexity of $I_A$ and $I_S$ and Jensen's inequality, we
have
\begin{equation*}
\int_0^1I_{A,S}(\chi_1'(s),\chi_2'(s))ds \leq
\int_0^1I_{A,S}(\phi_1'(s),\phi_2'(s))ds\:.
\end{equation*}
Since $\int_0^1I_{A,S}(\chi_1'(s),\chi_2'(s))ds=\tau
I_{A,S}(y)+(1-\tau)I_{A,S}(z)$ where $y=(y_1,y_2)$ and
$z=(z_1,z_2)$, it remains to minimize
\begin{equation*}
\delta q +\tau I_{A,S}(y)+(1-\tau)I_{A,S}(z)\:,
\end{equation*}
subject to
\[
\begin{cases} q\geq 0,\ \
y_1 \geq  y_2 \\
x_1=-q+(1-\tau)z_1+\tau y_2\\
x_2=q+(1-\tau)z_2\leq (1-\tau)z_1\\
q+(1-\tau)z_2 \leq (1-\tau)z_1. \\
\end{cases}
\]
Let $v_1=(1-\tau)z_1$ and $v_2=(1-\tau)z_2 +q$, we check that
$y_1=\frac{x_2-v_2}{\tau}$, $y_2=\frac{w+x_2-v_2}{\tau}$ and
\[J(x_1,x_2,w)=\inf\{\delta (w-x_1+x_2) +I_{A,S}(x_2,x_1), \inf_C
g(q,\tau,v_1,v_2)\}\:,\] $g$ et $C$ defined in Theorem
\ref{flle-action}.
\end{itemize}
\end{proof}

Let us give an interpretation of this result in a queueing
context. Indeed, the most likely ways to get departures with mean
$x_1$, times spent at the back of queue with mean $x_2$ and
waiting times with mean $w$ are when

\medskip

\begin{itemize}
\item all customers belong to the same busy period (the system is
never empty). Thus the mean of the departures is equal to the mean
of the services and the average time spent at the very back of the
queue is equal to the mean of the arrivals,

\medskip

\item or let $(q,\tau,v_1,v_2)\in C$ the variables where  $g$
reaches its minimum. Then the customer $-n\tau$ finds the system
empty, and the optimal path splits into two periods. Customer $-n$
waits $nq$ time before starting its service. During the first
period the average of arrivals and services are
$\frac{x_2-v_2}{\tau}$ and $\frac{x_2-v_2+w}{\tau}$ and it
consists of $n\tau$ customers. During the second period the mean
of services is $\frac{v_2-q}{1-\tau}$ and the arrivals occur with
mean $\frac{v_1}{1-\tau}$.
\end{itemize}

\medskip

The equation \eref{eq-foncinout} determines the rate functions of
the output variables in terms of the rate functions of the input
variables. A natural question is when are these functions equal?
More precisely, when do we have $(I_D,I_R)=(I_A,I_S)$? We start
with an example which translates the result stated in \cite{DMO04}
in terms of large deviations.

\medskip

\paragraph{Example: the M/M/1 queue} Suppose that
$\{a_n,\:n\in\Z\}$ and $\{s_n,\:n\in\Z\}$ are i.i.d. mutually
independent and for $t\geq 0$,
\[
P(a_0>t)=e^{-\lambda t}\:, \qquad P(s_0>t)=e^{-\mu t}\:.
\]
For $\theta<\lambda<\mu$,
$$\Lambda_A(\theta)=\log\Bigl(\frac{\lambda}{\lambda-\theta}\Bigr)\:, \qquad
\Lambda_S(\theta)=\log\Bigl(\frac{\mu}{\mu-\theta}\Bigr)\:.$$
Since
$\delta=\sup\{\theta>0,\:\Lambda_S(\theta)+\Lambda_A(-\theta)<0\}$,
we check that $\delta=\mu-\lambda$. Moreover we have
$$I_A(a)=\lambda a- \log(\lambda a)-1\:, \qquad I_S(s)=\mu s-
\log(\mu s)-1\:.$$ Resolving the above optimization problem, we
have
\begin{eqnarray*}
J(x_1,x_2,w)&=&\lambda x_1- \log(\lambda x_1)-1+\mu x_2- \log(\mu
x_2)-1+(\mu-\lambda)w\\
&=&I_{A,S}(x_1,x_2)+\delta w\:.
\end{eqnarray*}

    \subsection{Existence of fixed points}\label{sse-fix}
Suppose $\{a_n,\:n\in\Z\}$ and $\{s_n,\:n\in\Z\}$  i.i.d.,
mutually independent with means $1/\lambda$ and $1/\mu$. They
satisfy large deviations principles with rate functions $I_A$ and
$I_S$ differentiable on $\R_{+}$. Let $\beta>0$ such that
\begin{equation}\label{eq-Ia'=f(Is)'}
I_S'(x)-I_A'(x)=\beta,\qquad \forall x\in\R_+\:.
\end{equation}
We can rewrite this relation in the following which will be useful
for future analysis
\begin{equation}\label{eq-Ia=f(Is)}
I_A(x)-I_A(y)=I_S(x)-I_S(y)-\beta(x-y),\qquad \forall
(x,y)\in(\R_+)^2\:.
\end{equation}
Before stating the fixed point result for the Queue/Store operator
for the rate functions, we give an interpretation of the relations
\eref{eq-Ia'=f(Is)'} and \eref{eq-Ia=f(Is)} in terms of {\em
exponential tilting} of measures.

\medskip

\begin{defn}
Let $\beta>0$, we say that a measure $\nu'$ is the
$\beta$-exponentially-tilted measure of a given measure $\nu$ if
for all $x\in\R$,
\begin{equation}\label{eq-tilt}
\frac{d\nu'}{d\nu}(x)=\frac{e^{\beta x}}{\int e^{\beta x}
d\nu(x)}\:,
\end{equation}
where $(d\nu'/d\nu)(x)$ denotes the density of $\nu'$ with respect
to $\nu$.

\end{defn}

\medskip

\begin{prop}
Suppose $\{a_n,\:n\in\Z\}$ and  $\{s_n,\:n\in\Z\}$ satisfy
relation \eref{eq-Ia'=f(Is)'} or \eref{eq-Ia=f(Is)}, with
marginals $\nu_A$ and $\nu_S$. Then $\nu_A$ is the
$\beta$-exponentially-tilted function of $\nu_S$.
\end{prop}

\medskip

\begin{proof}
Recall that $I_A(1/\lambda)=0$. We check that
\begin{eqnarray}\nonumber
\Lambda_A(\theta)=\sup_{x\in\R}[\theta x-I_A(x)]
&=&\sup_{x\in\R}[\theta x-I_S(x)+I_S(1/\lambda)+\beta
(x-1/\lambda)]\\\nonumber &=&\sup_{x\in\R}[(\theta+\beta)
x-I_S(x)]-[\frac{\beta}{\lambda}-I_S(1/\lambda)]\\
\label{eq-AS} &=&\Lambda_S(\theta+\beta)-\Lambda_S(\beta) \:.
\end{eqnarray}
Considering the exponentials on both sides, we have
$$\int_{\R_+}e^{\theta x} d\nu_A(x)=\frac{\int_{\R_+}e^{(\theta+\beta) x}
d\nu_S(x)}{\int_{\R_+}e^{\beta x} d\nu_S(x)}\:.$$
\end{proof}

\medskip

\begin{lem}
Assume that $\{a_n,\:n\in\Z\}$ and $\{s_n,\:n\in\Z\}$ satisfy
assumptions $(i)$ and $(ii)$. Then the tilt coefficient $\beta$ is
given by
$$\beta=\sup\{\theta:\:\Lambda_S(\theta)+\Lambda_A(-\theta)\leq
0\}\:.$$
\end{lem}

\medskip

\begin{proof}
Let $\theta =-\beta$ in \eref{eq-AS}, we have
\[
\Lambda_A(-\beta)=\Lambda_S(0)-\Lambda_S(\beta)\:.
\]
Since $\Lambda_S(0)=0$, we obtain
$\Lambda_A(-\beta)+\Lambda_S(\beta)=0$. Thanks to \eref{eq-delta},
we have
$$\beta \leq\sup\{\theta\geq
0,\:\Lambda_S(\theta)+\Lambda_A(-\theta)\leq 0\}=\delta\:.$$ The
function $G(\theta)= \Lambda_S(\theta)+\Lambda_A(-\theta)$ is
convex with $G(0)=0$. Moreover, by assumption $(ii)$,
$$G'(0)=\Lambda'_S(0)-\Lambda'_A(0)< 0\:.$$ This means that $G$ is equal to zero twice
on $[0,\beta]$, at $0$ and $\beta$, then
$$\beta =\sup\{\theta\geq
0,\:\Lambda_S(\theta)+\Lambda_A(-\theta)\leq 0\}=\delta\:.$$
\end{proof}

\medskip

\begin{thm}\label{fix}
Let $\delta>0$ and the sequences $\{a_n,\:n\in\Z\}$ and
$\{s_n,\:n\in\Z\}$ the input variables of the model with $\nu_A$
the $\delta$-tilted measure of $\nu_S$ (or $\nu_S$ the
$(-\delta)$-tilted measure of $\nu_A$). Then the sequence
$\{(\frac{D_n}{n},\frac{R_n}{n},w_0/n),\:n\in\N^*\}$ satisfies a
large deviations principle on $\R_+^3$ with rate function
\begin{equation}
J(x_1,x_2,w)=\delta w+I_{A,S}(x_1,x_2)\:.
\end{equation}
Moreover, we have
$$I_D=I_A\quad \mrm{ and }\quad I_R=I_S\:.$$
\end{thm}

\medskip

\begin{proof}
By Theorem \ref{flle-action}, the sequence
$\{(\tilde{D}_n(1),\tilde{R}_n(1),w_1/n),\:n\in\N^*\}$ satisfies a
large deviations principle on $\R_+^3$ with rate function
\[
J(x_1,x_2,w)=\inf\{\delta (w-x_1+x_2)+I_{A,S}(x_2,x_1);\inf_C
g(q,\tau,v_1,v_2)\}\:,
\]
We start with the first term of this optimization problem, i.e.
$$x_1=\phi_2(1),\quad x_2=\phi_2(1)\quad q=w-x_1+x_2\:.$$
We have $w=x_1-x_2+q$, and thanks to \eref{eq-Ia=f(Is)},
\begin{eqnarray*}
\delta (w-x_1+x_2)+I_{A,S}(x_2,x_1)&=&\delta
(w-x_1+x_2)+I_{A,S}(x_1,x_2)+\delta (x_1-x_2)\\
&=&\delta w +I_{A,S}(x_1,x_2)\:.
\end{eqnarray*}
For the second term, let $(q,\tau,v_1,v_2)\in C$,
\begin{itemize}
\item[$\bullet$] If $v_1=v_2$ , from the convexity property, then
\begin{eqnarray*}
g(q,\tau,v_1,v_2)&=&\delta q +\tau
I_{A,S}\Bigl(\frac{x_2-v_2}{\tau},\frac{x_2-v_2+w}{\tau}\Bigr)
+(1-\tau)I_{A,S}\Bigl(\frac{v_2}{1-\tau},\frac{v_2-q}{1-\tau}\Bigr)\\
&\geq&\delta q +I_{A,S}(x_2,x_2+w-q) =\delta w
+I_{A,S}(x_1,x_2)\:.
\end{eqnarray*}
\item[$\bullet$] If $v_1>v_2$, then
$$I_{A,S}\Bigl(\frac{x_2-v_2}{\tau},\frac{x_2-v_2+w}{\tau}\Bigr)=
I_{A,S}\Bigl(\frac{x_2-v_2+w}{\tau},\frac{x_2-v_2}{\tau}\Bigr)
+\delta\frac{w}{\tau}\:$$ Notice that $$
I_S\Bigl(\frac{v_2-q}{1-\tau}\Bigr)=
I_S\Bigl(\frac{v_2}{1-\tau}\Bigr) +
I_A\Bigl(\frac{v_2-q}{1-\tau}\Bigr)
-I_A\Bigl(\frac{v_2}{1-\tau}\Bigr) -\delta\frac{q}{1-\tau}.$$
Since $I_A$ is convex and $v_1 \geq v_2$, we have
$$I_A\Bigl(\frac{v_2-q}{1-\tau}\Bigr)-I_A\Bigl(\frac{v_2}{1-\tau}\Bigr)
\geq
I_A\Bigl(\frac{v_1-q}{1-\tau}\Bigr)-I_A\Bigl(\frac{v_1}{1-\tau}\Bigr).$$
Thus,
$$I_S\Bigl(\frac{v_2-q}{1-\tau}\Bigr)+I_A\Bigl(\frac{v_1}{1-\tau}\Bigr)
 \geq I_S\Bigl(\frac{v_2}{1-\tau}\Bigr)
+I_A\Bigl(\frac{v_1-q}{1-\tau}\Bigr) -\delta \frac{q}{1-\tau}.$$
Since $I_A$ and $I_S$ are convex and $x_2-v_2+w+v_1-q=x_1$, we
have
\begin{equation*}
g(q)\geq\delta w +I_{A,S}(x_1,x_2)\:.
\end{equation*}
We conclude that $J(x_1,x_2,w)=\delta w +I_{A,S}(x_1,x_2)$.
\end{itemize}
\end{proof}

In other words, the above theorem states that if the input
variables satisfy large deviations principles with rate functions
$I_A$ and $I_S$ with $\nu_A$ the $\delta$-tilted measure of
$\nu_S$, then the output variables satisfy large deviations
principles with the same rate functions.

\end{document}